\documentclass[11pt]{article}
\usepackage{amsmath, amsthm}
\usepackage{url}
\usepackage{amsfonts,amssymb}
\usepackage{epsfig}
\usepackage{graphicx}
\newcommand{\bm}{\mathbf}
\newcommand{\B}{\bigskip}
\newcommand{\m}{\medskip}
\newtheorem{theorem}{Theorem}
\newtheorem{cor}{Corollary}
\newtheorem{lemma}{Lemma}

\newtheorem{prop}{Proposition}
\theoremstyle{definition}
\newtheorem{definition}[theorem]{Definition}
\newtheorem{example}[theorem]{Example}

\numberwithin{equation}{section}

\date{ }

\title{The Number of Ways to Assemble a Graph} 
\author{Andrew Vince and Mikl\'os B\'ona\\ University of Florida, Department of Mathematics 
\\ Gainesville, FL, USA \\ 
  {\tt  avince@ufl.edu, bona@ufl.edu} 
}

\begin{document}  
\maketitle
 
\begin{abstract}  Motivated by the question of how macromolecules assemble, the notion of an {\it assembly tree}
of a graph is introduced.  Given a graph $G$, the paper is concerned with enumerating the number of assembly trees of $G$, a problem that applies to the macromolecular assembly problem.   Explicit formulas or generating functions are provided for the number of assembly trees of several families of  graphs, in particular for what we call $(H,\phi)$-graphs.  
In some natural special cases, we  apply powerful recent results of Zeilberger and Apagodu on multivariate generating functions, and results of Wimp and Zeilberger, to deduce recurrence relations and very precise asymptotic formulas for the number of assembly trees of the complete bipartite graphs $K_{n,n}$ and the complete tripartite graphs $K_{n,n,n}$.
Future directions for reseach, as well as  open questions, are suggested.
\end{abstract} 



\section{Introduction} \label{secIntro}

Although the context of this paper is graph theory, the concept of an assembly tree originated in an attempt to 
understand macromolecular assembly \cite{BSV}.  The capsid of a virus - the shell that protects the genomic material - self-assembles spontaneously, rapidly and quite accurately in the host cell.  Although the structure of the capsid is fairly well known, the assembly process by which hundreds of subunits (monomers) interact to form the capsid is not well understood.  In many cases, the capsid can be modeled by a polyhedron, the facets representing the monomers. The assembly of the capsid can be modeled by a rooted tree, the leaves representing the facets, the root the completed polyhedron, and the internal nodes intermediate subassemblies. The enumeration of such trees plays a central role in understanding how
symmetry effects  the assembly process \cite{BSV}. \m

All graphs in this paper are simple.  Let $G=(V,E)$ be a connected graph of order $n$ with vertex set $V$ and edge set $E$.   In the definition of an
assembly tree $T$ for the graph $G$, each node of $T$ is labeled by a subset of $V$. No distinction will be made between the node and its label. For a node $U$ in a rooted tree,  $c(U)$ denotes the set of children of $U$.

\begin{definition} An {\em assembly tree} for a connected graph $G$ on $n$ vertices is a rooted tree, each node of which is labeled by a subset $U \subset V$ such that 
\begin{enumerate}
\item each internal (non-leaf) node has at least two children,
\item there are $n$ leaves with labels  $\{v\},  \, v \in V$,
\item the label on the root is $V$,
\item  $U = \bigcup c(U)$ for for each internal node $U$.   
\end{enumerate} 
\end{definition}

\noindent An assembly tree $T$ for $G$ describes a process by which $G$ assembles.  At the 
beginning are the individual vertices of $G$ - the leaves of $T$.  Each internal node $U$ of $T$
represents the subgraph of $G$ induced by the subset $U$ of vertices.  Each internal node
also represents the stage in the assembly process by which  subgraphs of $G$ join
to form a larger  subgraph; more precisely, the subgraphs induced by the children of $U$ join to form the subgraph
induced by $U$. The process  terminates at the root of $T$ - representing the entire graph $G$.  Call two assembly trees $T_1$ and $T_2$ for a graph $G$ {\em equal} if there is a label  preserving graph isomorphism between $T_1$ and $T_2$.  \m

There are numerous ways, some mentioned in the last section, to further restrict how the assembly process occurs.
In this paper we will assume that, at each stage, two subgraphs can be joined if and only if there is an edge that connects them. 

\begin{definition}  \label{def:edge} An assembly tree for a connected graph $G = (V,E)$ using the {\em edge gluing rule}  is an assembly tree for $G$ satisfying the additional property:
\begin{enumerate}
\item[5.] Each internal node has exactly two children, and
if $U_1$ and $U_2$ are the
children of internal node $U$, then there is an edge $\{v_1,v_2\} \in E$, the {\em gluing edge}, 
such that $v_1 \in U_1$ and $v_2 \in U_2$.
\end{enumerate} \m
 \end{definition}

\begin{figure}[htb] \label{fig1}
\vskip -.5cm
\begin{center}
\includegraphics[width=10cm, keepaspectratio]{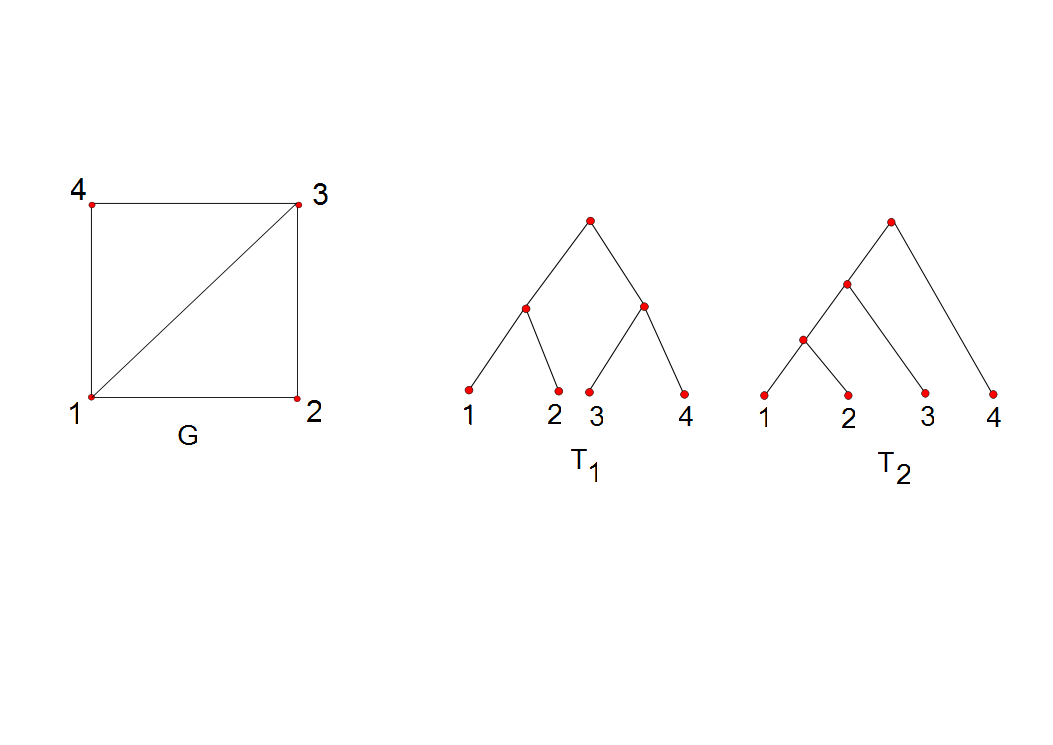}
\vskip -2.5cm
\caption{Two assembly trees for the graph $G$.}
\end{center} 
\end{figure}

\noindent  Figure~\ref{fig1} shows a graph $G$ and two assembly trees for $G$ using the edge gluing rule.  Throughout this paper, until the last section, all assembly trees use the edge gluing rule.  Therefore, the term ``assembly tree" will refer to an assembly tree using the edge gluing rule. 

The subject of the paper is, given a graph $G$,  to enumerate the number $$a(G)$$ of assembly trees of $G$. Gluing sequences are defined in  Section~\ref{secPathCycle} and are used too enumerate the number of assembly trees for paths, cycles and certain star graphs.  The concept of an $H$-graph is defined in Section~\ref{secH-partite}; the complete multi-partite graphs are special cases.   A generating function formula for the number of assembly trees for any $H$-graph is  provided in  Section~\ref{secH-partite}.  Section \ref{secExamples} considers three specific examples of $H$-graphs which
lead to frequently encountered families of graphs, such as complete bipartite graphs or complete tripartite graphs. 
For each of these examples, the relevant multivariate generating function is computed, then the {\em diagonal} of that
generating function is introduced and 
 studied. Very strong recent results of Doron Zeilberger and Moa Apagodu enable us to prove polynomial recurrence
relations for the coefficients of these diagonals, while results of Zeilberger and Jet Wimp allow us find the growth
rate of these coefficients at an arbitrary level of precision. In particular, we obtain the growth rates for the number of assembly trees of the complete bipartite graphs $K_{n,n}$ and the complete tripartite graphs $K_{n,n,n}$ as a function of $n$.  Open questions and further research directions are offered in Section~\ref{secOpen}.

\section{Paths, Cycles and Stars} \label{secPathCycle}

An assembly tree for a graph $G = (V,E)$ of order $n$, as defined in the introduction,  is a binary tree with $n$ leaves and $n-1$ internal nodes.  To each internal node $U$ there is a corresponding gluing edge as in Definition~\ref{def:edge}  (not necessarily unique), which we denote by  $e_U \in E$.

\begin{lemma} If $T$ is an assembly tree for a connected graph $G = (V,E)$, then the set of gluing edges
$\{e_U \, | \, U \text{  is an internal node of assembly tree  } T\}$ is a spanning tree of $G$.  
\end{lemma}

\proof If the set $S := \{e_U \, | \, U \text{  is an internal node of assembly tree  } T$\} of gluing edges is not spanning, then the root of $T$ would not be $V$.  If $S$ contains a cycle, then $T$ would have a node with just one child.  If $S$ is not connected, then $G$ would not be connected.
\qed \m

If $S \subseteq E(G)$ is any spanning tree of a connected graph $G$, then any linear ordering $e_1, e_2, \dots , e_{n-1}$ of the edges in $S$ induces an assembly tree for $G$ as 
follows.  Build the tree $T$ from the bottom up.  The leaves are the singleton vertices of $G$.  Assume that we have proceeded through the sequences of edges from $e_1$ to $e_{k-1}$.  For $e_k = \{u_1,u_2\}$ add a node to $T$ whose two children are the already constructed nodes $U_1$ and $U_2$ such that $u_1 \in U_1$
and $u_2 \in U_2$.  Call an ordering  $e_1, e_2, \dots , e_{n-1}$ of the edges of the spanning tree a {\em gluing sequence}.  The elements of a gluing sequence for $G$ are the gluing edges of the corresponding assembly tree $T$.  

\begin{figure}[htb] \label{fig2}
\vskip -1cm
\begin{center}
\includegraphics[width=4cm, keepaspectratio]{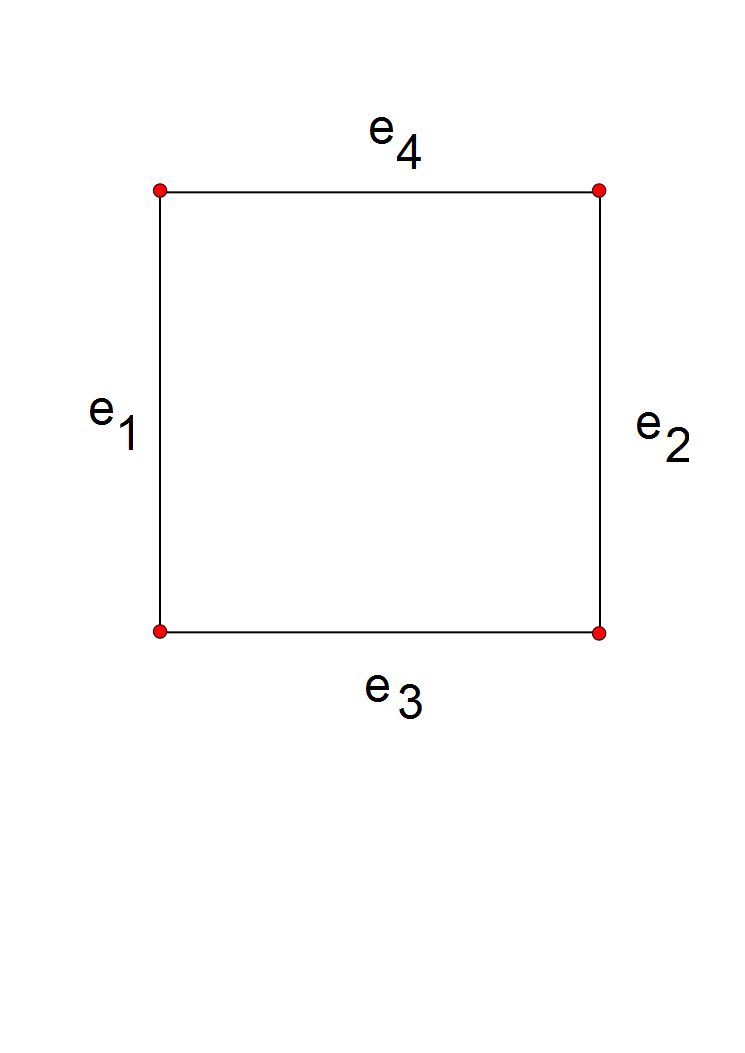}
\vskip -1.5cm
\caption{The cycle $C_4$.}
\end{center} 
\end{figure}

\begin{example} Consider the $4$-cycle $C_4$ in Figure~\ref{fig2}.  This example shows that 
 two different spanning trees can induce the same assembly tree: the gluing sequences  $(e_1,e_2,e_3)$ and $(e_1,e_2, e_4)$ produce the same assembly tree. Moreover, two different orderings, for example $(e_1,e_2,e_3)$ and $(e_2,e_1,e_3)$, of the same spanning tree can produce the same assembly tree. \end{example} 

For the star $S_n := K_{1,n}$, the spanning tree is $S_n$ itself, and each gluing sequence produces a distinct assembly tree.  This leads immediately to the following result.

\begin{prop} For the star the number of assembly trees is  $a(S_n) = n!$.
\end{prop}

Consider the star $S^2_n$ with $n$ arms such that each arm has length $2$, as in Figure 3.  

\begin{figure}[b] \label{fig3}
\begin{center}
\includegraphics[width=7cm, keepaspectratio]{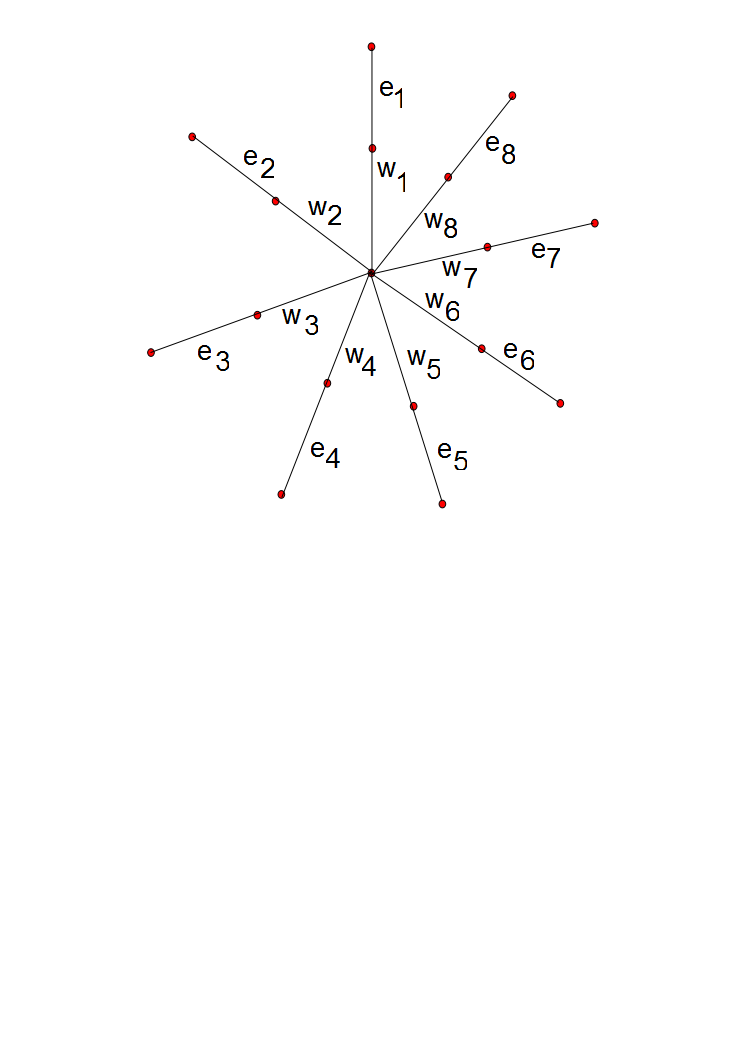}
\vskip -5cm
\caption{The star $S^2_8$.}
\end{center} 
\end{figure}

\begin{theorem} $$a(S^2_n) = \sum_{k=0}^n \binom{n}{k} \,\frac{(2n-k)!}{2^{n-k}}$$
\end{theorem}

\proof  Suppose that $k$ of the $e$-edges come first  in the gluing sequence.  In Figure~\ref{fig3} we refer to the $e$-edges and the $w$-edges.  There are $\binom{n}{k}$ ways to 
choose these edges, and the order does not matter for the assembly tree.  There are $2n-k$ edges that remain.
For convenience label them $e_1,e_2, \dots, e_{n-k}$ and $w_1, w_2, \dots , w_n$.  These $2n-k$ edges
can be placed in any order in the gluing sequence as long as the edge $w_i$ comes after $e_i$ for $i = 1,2, \dots, n-k$.   Each such order determines a distinct assembly tree.  To determine the number of such permutations, first choose $k$ positions for $w_{n-k+1}, \dots , w_n$.  There are $\binom{2n-k}{k}$ ways to do this, and for each such choice 
the edges $w_{n-k+1}, \dots , w_n$ can be permuted in $k!$ ways.  The remaining $2(n-k)$ positions are to be filled by the edges $e_1,e_2, \dots, e_{n-k}$ and  $w_1, w_2, \dots , w_n$ so that  the edge $w_i$ comes after $e_i$ for $i = 1,2, \dots, n-k$.  The number of ways to do this equals the number of permutations of $2(n-k)$ objects where there are
$2$ objects of type $1$, $2$ objects of type $2\dots $, $2$ objects of type $n-k$.  This is equal to 
$\frac{(2n-2k)!}{2^{n-k}}$.  So, with  $k$ of the $e$-edges coming first in the gluing sequence there are
$$ \binom{n}{k} \,  \binom{2n-k}{k} \, k! \frac{(2n-k)!}{2^{n-k}} = \binom{n}{k} \,\frac{(2n-k)!}{2^{n-k}}$$
assembly trees.  Summing over all possible values of $k$ from $k=0$ to $k=n$ gives the formula in the statement 
of the theorem.
\qed \m 

\begin{theorem} If $P_n$ is the path and $C_n$ is the cycle on $n$ vertices, then
$$a(P_n) = \frac{1}{n} \binom{2n-2}{n-1}, \qquad \qquad a(C_n) =  \frac12 \,\binom{2n-2}{n-1}.$$ 
\end{theorem}

\proof For the path, the unique spanning tree $S$ consists of all the edges of the path.  We proceed 
by induction.  First consider the number of assembly trees in the case that $e \in S$ is the last edge in the gluing sequence. If the removal of $e$ from $G$ results in subgraphs of orders $k$ and $n-k$, then the number of assembly trees such that $e$ is the last edge in the gluing sequence is $a(P_k) \, a(P_{n-k})$.  If $T$ and $T'$ are assembly trees coming from gluing sequences with distinct last elements, then $T \neq T'$.  Therefore $a(P_n) = \sum_{k=1}^{n-1} \, a(P_k) \, a(P_{n-k})$, which is a well know recurrence for the Catalan numbers. 

Concerning the cycle, there are $n$ spanning trees of $C_n$. Given any one of these spanning 
trees, by the result above for the path, there are $ \frac1n \,\binom{2n-2}{n-1}$ corresponding assembly trees for
$C_n$, hence a total of $n \, \frac{1}{n} \binom{2n-2}{n-1} = \binom{2n-2}{n-1}$  assembly trees.  But each of these is counted twice for following reason.  The assembly tree corresponding to a sequence of edges in a spanning tree
for which $e$ is the last edge in the gluing sequence and $f$ is the edge of $G$ not in the gluing sequence is equal to the
assembly tree for which $f$ is the last edge in the gluing sequence and $e$ is the edge of $G$ not in the gluing sequence.
\qed

\section{$H$-graphs} \label{secH-partite}

Let $H$ be a connected graph with vertex set $[N] := \{1,2,\dots , N\}$, and let $\phi \, : \, [N] \rightarrow \{0,1\}$ be a  labeling of the vertices of $H$.  For any sequence $(n_1, n_2, \dots, n_N)$  of non-negative integers, define a graph
$G_{(H,\phi)}(n_1, n_2, \dots, n_N)$ as follows.  The vertex set is 
$$V(G_{(H,\phi)}(n_1, n_2, \dots, n_N)) :=  \{ (i,j) \,: \,  i \in [N], \,  1\leq j \leq n_i\}$$
 and $(i,j)$ is adjacent to $(i',j')$ if and only if 
$$\begin{aligned} i &= i' \quad \text{and} \quad \phi(i) = 1, \quad \text{or} \\
i &\neq i' \quad \text{and} \quad \{i,i'\} \in E(H).
\end{aligned}$$

\noindent This is equivalent to saying that the graph $G_{(H,\phi)}(n_1, n_2, \dots, n_N)$  is obtained by replacing each vertex $i$ of $H$ by a complete graph of order $n_i$ or its complement  ($n_i$ isolated vertices), and by  replacing each edge of $H$ by all possible edges  between the graphs that replace the two end vertices of that edge.     Call  a graph  an $ (\bm H,\phi)-${\em graph} if it is of the form $G_{(H,\phi)}(n_1, n_2, \dots, n_N)$ for some choice of the parameters $\{n_1,\dots, n_N\}$. \m

The following notation will be used in this section
\begin{enumerate} 
\item${\bm n} := (n_1,n_2, \dots , n_N)$  
\item ${\bm n} \geq {\bm k} \quad \text{if and only if} \; n_i \geq k_i \; \text{for all}\; i$ 
\item ${\bm x} := (x_1, x_2, \dots , x_N)$ 
\item ${\bm n}!  = n_1!n_2!\cdots n_N!$ 
\item $\binom{\bm n}{\bm k} = \binom{n_1}{k_1} \binom{n_2}{k_2} \cdots \binom{n_N}{k_N}$ 
\item ${\bm x}^{\bm n} = x_1^{n_1} x_2^{n_2} \cdots x_N^{N^n}$
\item $a_{(H,\phi)} ({\bm n}) = a(G_{(H,\phi)}(n_1, n_2, \dots, n_N))$ 
\item $A_{(H,\phi)}(\bm x) = \sum_{{\bm n} \geq {\bm 0}} \,  a_{(H,\phi)} ({\bm n}) \, \frac{{\bm x}^{\bm n}}{{\bm n}!}$.
\end{enumerate}
\noindent  The last entry in the above list is the exponential  generating function for the number of assembly trees
of an $(H,\phi)$-graph. The zero vector is denoted $\bm 0$.    

\begin{theorem}  \label{thmMain}
The exponential generating function for a connected $(H,\phi)$-graph is
 $$A_{(H,\phi)}(\bm x) = 1 - \sqrt{1-2\sum_{i=1}^{N}  x_i+ \sum_{\phi(i) = 0} x_i^2 +2 \sum_{\{i,j\} \notin E(H)} x_i x_j }.$$ \end{theorem}

\proof Let $\bm 0$ be the all zeros vector;  let ${\bm e_i}$ be the vector with each coordinate $0$ except the $i^{th}$ coordinate $1$; and let $\bm e_{i,j}$ be the vector with each coordinate $0$ except the $i^{th}$ and $j^{th}$ coordinate $1$.  Note
that 
$$\begin{aligned} a_{(H,\phi)}(\bm 0) &= 0 \\
a_{(H,\phi)}({\bm e_i}) &= 1 \quad \text{for all} \; i \\
a_{(H,\phi)}(2\,{\bm e_i}) &= 0 \quad \text{if}\; \phi(i) = 0 \\
a_{(H,\phi)}({\bm e_{i,j}}) & = \begin{cases} 1 \quad \text{if} \; \{i,j\}\in E(H) \\ 0 \quad \text{if} \; \{i,j\}\notin E(H)  \end{cases} 
\end{aligned}$$
The following recurrence holds for all $\bm n$ except those of the form $\bm e_i, \, \bm e_{i,j}$ when $\{i,j\} \notin E(H)$ and $2\bm e_i$ when $\phi(i) = 0$.  To simplify notation, denote  $a_{(H,\phi)} ({\bm n})$ by $a({\bm n})$ and
$ A_{(H,\phi)}(\bm x)$ by  $A(\bm x)$. 
$$a({\bm n}) = \frac12\, \sum_{\bm 0 \leq \bm k \leq \bm n} \binom{{\bm n}}{{\bm k}} a({\bm k})a({\bm n-\bm k}).$$
The recurrence above is obtained by considering the two subtrees $T_1$ and $T_2$, rooted at each of the children of the root of an assembly tree. The tree $T_1$ is itself an assembly tree of a graph of the form $G_{(H,\phi)}(k_1, k_2, \dots, k_N)$ where $k_i \leq n_i$ for all $i$, and $T_2$ is an assembly tree of a graph of the form 
 $G_{(H,\phi)}(n_1-k_1, n_2-k_2, \dots, n_N-k_N)$.  
Now
$$\begin{aligned}  A(\bm x) &= \sum_{{\bm n} \geq {\bm 0}} \,  a({\bm n}) \, \frac{{\bm x}^{\bm n}}{{\bm n}!} \\
&=  \sum_{{\bm n} \geq {\bm 0}} \,  \left \{ \frac12\, \sum_{\bm 0 \leq \bm k \leq \bm n} \binom{{\bm n}}{{\bm k}} a({\bm k})a({\bm n-\bm k}) \right \} \, \frac{{\bm x}^{\bm n}}{{\bm n}!}  +\sum_{i=1}^{N}  x_i - \sum_{\phi(i) = 0} \frac{x_i^2}{2} - \sum_{\{i,j\} \notin E(H)} x_i x_j  \\
&= \frac12\, A^2(\bm x) +\sum_{i=1}^{N}  x_i - \sum_{\phi(i) = 0} \frac{x_i^2}{2} - \sum_{\{i,j\} \notin E(H)} x_i x_j. 
\end{aligned}$$
The added terms  $\sum_{i=1}^{N}  x_i - \sum_{\phi(i) = 0} \frac{x_i^2}{2} - \sum_{\{i,j\} \notin E(H)} x_i x_j$ correct for the three cases for which the recurrence does not hold.  Solving for $A(\bm n)$ by the quadratic formula yields the
generating function in the statement of the theorem.  
\qed \m 

Some special cases follow as corollaries. For example, if $H$ is just a single vertex $v$ and $\phi(v) =1$, then $G_{(H,\phi)}(n)$ is the complete graph.   Therefore, by Theorem~\ref{thmMain} the generating function for the complete graph is $1-\sqrt{1-2x}$, which when expanded gives the following result. 

\begin{cor} If $K_n$ is the complete graph, then
$$a(K_n) = \frac{(2n-2)!}{2^{n-1}\, (n-1)!}.$$
\end{cor}

\begin{cor}  \label{multipartite}
 If $K_{(n_1, n_2, \dots, n_N)}$ is the complete multipartite graph, then its exponential generating 
function is 
$$A(\bm x) = 1 - \sqrt{1-N + \sum_{i=1}^{N} \, (1-x_i)^2}.$$
In particular, the generating function for the number of assembly trees of the complete bipartite graph $K_{n_1,n_2}$ is 
$$A(x,y) = 1-\sqrt{(1-x)^2 + (1-y)^2 - 1}.$$
\end{cor} 

\proof If $H$ is $K_n$ and $\phi$ is identically $0$ then $G_{(H,\phi)}(n_1, n_2, \dots, n_N)$ is the complete bipartite graph. \qed \m

Expanding for the generating function $A(x,y)$ counting assembly trees of 
 complete bipartite graphs using Maple we arrive at the  exponential generating function

$$\begin{aligned} A(x,y) = & \; y+x+xy+xy^2+x^2y+xy^3+x^3y+(5/2)x^2y^2+(9/2)x^2y^3  
                                      +x^4y  \\ & +(9/2)x^3y^2+ xy^4+7x^2y^4+7x^4y^2+(25/2)x^3y^3+10x^2y^5  
                                      +(55/2)x^4y^3 \\ & +(55/2)x^3y^4 +10x^5y^2+(27/2)x^2y^6+(645/8)x^4y^4+(105/2)x^3y^5  \\ & 
                                      +(27/2)x^6y^2+(105/2)x^5y^3 +yx^5+y^5x+yx^6+y^6x+yx^7+y^7x+yx^8   \\ & +(35/2)y^2x^7 
+91y^3x^6 +(1575/8)y^4x^5 +(1575/8)y^5x^4  +91y^6x^3 \\ & +(35/2)y^7x^2+y^8x \cdots, \end{aligned}$$
which gives the following table for the  number of assembly trees.  
  
$$\begin{matrix}  1 & 2 & 6 & 24 \\  & 10 & 54 & 336\\& & 450 & 3960 \\ & & & 46400  \end{matrix}$$

The diagonal elements $1,10,450, 23200, \dots$ are the number of assembly trees for $K_{1,1}, K_{2,2}, K_{3,3}, K_{4,4}, \dots$, a sequence which does not match anything in the Online Encyclopedia of Integer Sequences \cite{OEIS}.  We will return to this sequence in the next section, where we will
 find the asymptotic growth rate of the sequence, and to find a polynomial recurrence
relation satisfied by the sequence.

Consider the set of all $(H,\phi)$-graphs on $n$ labeled vertices.  In other words, such a graph is obtained
by choosing, say $n_1$ of the $n$ vertices to correspond to vertex $1$ of $H$,  $n_2$ of the $n$ vertices to correspond to vertex $2$ of $H$, ... ,  $n_N$ of the $n$ vertices to correspond to vertex $N$ of $H$.  Let $b_{(H,\phi)}(n)$ denote the total number of assembly trees of all the  possible $(H,\phi)$-graphs on $n$ labeled vertices, and $B_{(H,\phi)}(n)$ the corresponding exponential generating function
$$B_{(H,\phi)}(n) := \sum_{n=0}^{\infty}  b_{(H,\phi)}(n) \, \frac{x^n}{n!}.$$

\begin{cor} The exponential generating function of the number or assembly trees of all connected $(H, \phi)$-graphs 
of order $n$ such that the number of vertices in $H$ is $N$, the number of edges in $H$ is $M$, and the
number of $i$ such that $\phi(i)=0$ is $J$ is
$$B_{(H,\phi)}(n) =  1 - \sqrt{1-2N x+ \left ( 2\binom{N}{2} - 2M +J \right ) x^2}.$$ 
\end{cor}

\proof We have
$$\begin{aligned} 
 \sum_{i=0}^{\infty}  b_{(H,\phi)}(n) \, \frac{x^n}{n!} &= \sum_{n=0}^{\infty} \left [ \sum_{n_1+n_1+\cdots+n_N = n} \binom{n}{n_1 n_2 \cdots  n_N}  a_{(H,\phi)} ({\bm n}) \right ] \, \frac{x^n}{n!} \\
&= \sum_{n=0}^{\infty} \left [ \sum_{n_1+n_1+\cdots+n_N = n} \frac{ a_{(H,\phi)} ({\bm n})}{n_1! n_2!\cdots n_N!}
\right ]\,  x^n \\ &=  \sum_{{\bm n} \geq {\bm 0}} \,  a_{(H,\phi)} ({\bm n}) \, \frac{x^ n}{{\bm n}!} \\
&=  \sum_{{\bm n} \geq {\bm 0}} \,  a_{(H,\phi)} ({\bm n}) \, \frac{x^ {n_1 + \cdots +n_N}}{{\bm n}!} \\
&= A_{(H,\phi)}(x,x,\dots , x) =  1 - \sqrt{1-2\sum_{i=1}^{N}  x+ \sum_{\phi(i) = 0} x^2 +2 \sum_{\{i,j\} \notin E(H)} x^2} \\
&=  1 - \sqrt{1-2N x+ \left ( 2\binom{N}{2} - 2M +J \right ) x^2},
\end{aligned}$$ 
where the second-to-last equality follows from Theorem~\ref{thmMain}. 
\qed

\section{Examples} \label{secExamples}
In this section, we consider a few interesting examples. In these examples, the graph $H$ is the basis for the
construction of an $H$-graph will be very small (two or three
vertices), and $n_1=n_2$ or $n_1=n_2=n_3$ will hold, resulting in graphs $G_{(H,\phi)}$ with two or three classes
of vertices in some obvious sense. 

\subsection{Theoretical Background}
\subsubsection{Power Series in One Variable}

Let $\mathbb{C}[n]$ denote the ring of all polynomials in one variable over the field of complex numbers, and
let $\mathbb{C}[[x]]$ denote the ring of all formal power series with complex coefficients. 
In what follows, we present a few important definitions and theorems on one-variable power series. 
The interested reader can consult Chapter 6 of \cite{StanleyEC2} for a deeper introduction to the topic,
including the proofs of the theorems we include here.  

\begin{definition} \label{prekdef}
 A  sequence $f(0),f(1),\cdots $ of complex numbers is called {\em polynomially recursive}, or {\em $p$-recursive}
if there exist polynomials $P_0,P_1,\cdots, P_k\in \mathbb{C} [n]$, with $P_k\neq
0$ so that \begin{equation} \label{rekudef}
P_k(n+k)f(n+k)+P_{k-1}(n+k-1)f(n+k-1)+\cdots +P_0(n)f(n)=0
\end{equation} for all 
natural numbers $n$. 
\end{definition}

\begin{definition} We say that the power series $u(x)\in \mathbb{C} [[x]] $
is $d$-finite if there
exists a positive integer $d$ and polynomials  $p_0(n),p_1(n),\cdots ,
p_d(n)$  so that $p_d\neq 0$ and
\begin{equation} \label{findef} p_d(x)u^{(d)}(x)+p_{d-1}(x)u^{(d-1)}(x)+\cdots
+p_1(x)u'(x)+p_0(x)u(x)=0, \end{equation} 
Here $u^{(j)}=\frac{d^ju}{dx^j}$.
\end{definition}

\begin{theorem} \label{dfinprec}
 The sequence $f(0), f(1), \cdots $ is $p$-recursive if and
only if its ordinary generating function 
\begin{equation} 
\label{dgenfun} u(x)=\sum_{n=0}^\infty f(n)x^n \end{equation}
is $d$-finite.
\end{theorem}

\begin{definition} \label{algdef}
The formal power series $f\in \mathbb{C}[[x]]$ is called
{\em algebraic} if there exist polynomials $P_0(x),P_1(x),\cdots ,P_d(x)
\in \mathbb{C}[x]$ that are not all equal to zero so that
\begin{equation} \label{eqalgdef} 
P_0(x)+P_1(x)f(x)+\cdots +P_d(x)f^d(x)=0.
\end{equation}
The smallest positive $d$ for which such polynomials exist is called the
{\em degree} of $f$.
\end{definition}

\begin{theorem} \label{algdif} If $f\in \mathbb{C}[[x]]$
 is algebraic, then it is $d$-finite.
\end{theorem}

 We point out that the converse of 
Theorem \ref{algdif} is not true. For instance, $f(x)=\sum_{n\geq 1} \frac{x^n}{n}=\ln(1/(1-x))$ is
$d$-finite, but not algebraic, as we will soon see. 

One way to prove that a power series is not algebraic is by proving that it is not $d$-finite. Another way of proving that a power series is  not algebraic is by showing that is does not have the
"right" growth rate. The following theorem of Jungen is a powerful tool in doing so. 

\begin{theorem} \cite{Jungen} \label{jungen}
Let $f(x)=\sum_{n\geq 0}a_nx^n\in \mathbb{C}[[x]]$ be an {\em algebraic} power series,
 and let us assume that $a_n\sim c n^r \alpha^n$, where $c$ and
$\alpha$ are {\em non-zero}
complex constants, and $r$ is a {\em negative real} constant.  

Then $r=s+\frac{1}{2}$, for some negative integer $s$. 
\end{theorem}
In particular, selecting $c=1$, $r=-1$ and $\alpha=1$, we see that $f(x)=\sum_{n\geq 1} \frac{x^n}{n}$ is not algebraic.

\subsubsection{Power Series in Several Variables}
Now we consider formal power series in several variables. For a deeper introduction to the topic, including the
proofs of the theorems we present, see \cite{lipshitz}. 
Let $\mathbb{C}[[x_1,x_2,\cdots ,x_k]]$ denote the algebra of all formal power series in variables
$x_1,x_2,\cdots ,x_k$ over the field of complex numbers. 
\begin{definition}
Let $f(n_1,n_2,\cdots ,n_k): \mathbb{N} ^k \rightarrow \mathbb{C}$ be a function, 
and let $F(x_1,x_2,\cdots ,x_k)=\sum_{n_1,n_2,\cdots ,n_k}f(n_1,n_2,\cdots ,n_k)x_1^{n_1}x_2^{n_2}\cdots 
x_k^{n_k}\in \mathbb{C}[[x_1,x_2,\cdots ,x_k]]$. 
\vskip 1mm
We say that $F$ is $d$-finite if all the derivatives
\[ \left(\frac{\partial }{\partial x_1}\right)^{d_1} \left(\frac{\partial }{\partial x_2}\right)^{d_2} \cdots
\left(\frac{\partial }{\partial x_k}\right)^{d_k}F\]
for $d_i\geq 0$ lay in a finite dimensional vector space over the field of rational functions $\mathbb{C}(x_1,x_2,\cdots
,x_k)$.
\end{definition}

\begin{theorem} \label{algtodfinmult}
 Let $F\in \mathbb{C}[[x_1,x_2,\cdots ,x_k]]$. If $F$ is algebraic, then it is $d$-finite.
\end{theorem}

The notion of the {\em diagonal} of a multivariate power series $F$ is a natural one in that it enables us to 
focus on the coefficients of $F$ that are often the most interesting for practical purposes. 

\begin{definition}
Let $f(n_1,n_2,\cdots ,n_k): \mathbb{N} ^k \rightarrow \mathbb{C}$ be a function, 
and let $F(x_1,x_2,\cdots ,x_k)=\sum_{n_1,n_2,\cdots ,n_k}f(n_1,n_2,\cdots ,n_k)x_1^{n_1}x_2^{n_2}\cdots 
x_k^{n_k}$. Then the {\em diagonal} of the multivariate power series $F(x_1,x_2,\cdots ,x_k)$ is the
{\em univariate} power series
\[\rm{diag} F(x)=\sum_n f(n,n,\cdots ,n) x^n.\]
\end{definition}

\begin{example} Let $F(s,t)=\frac{1}{1-s-t}=\sum_{m\geq 0} (s+t)^m$. Then for every $n\in \mathbb{N}$, the coefficient
of $s^nt^n$ in $F(s,t)$ is equal to ${2n\choose n}$. Therefore, 
\[\rm{diag} F(x) = \sum_n {2n\choose n} x^n = \frac{1}{\sqrt{1-4x}}.\]
\end{example}

\begin{theorem}  \label{diagdfin}
Let $F\in \mathbb{C}[[x_1,x_2,\cdots ,x_k]]$. If $F$ is $d$-finite, then $\rm{diag}F(x)$ is
also $d$-finite.
\end{theorem}

\subsection{Three families of graphs}

Our first example of computing the  diagonal of a power series $A_{(H,\phi)}(\bf{z})$ is included because of
the precise nature of the answer that we are able to compute. 

\begin{example} \label{easyex}
 Let $H$ be a graph on vertex set $\{u,v\}$, with $\phi(u)=0$ and $\phi(v)=1$, and with one edge,
the edge $uv$. The graphs $G_{(H,\phi)}$ are the graphs consisting of a subgraph $G'$ consisting of $n_1$ independent vertices and a complete subgraph $G''$ on $n_2$ vertices so that $G'$ and $G''$ are vertex-disjoint, and  each
vertex of $G'$ is adjacent to each vertex of $G''$. 
Then by Theorem \ref{thmMain}, we have
\[A_{(H,\phi)}(x,y)=1-\sqrt{1-2x-2y+y^2}.\]
\end{example}

In particular, the diagonal of $A_{(H,\phi)}(x,y)$ counts the number of assembly trees of such 
graphs with $n_1=n_2$. In order to compute this diagonal, note that
\[A_{(H,\phi)}(x,y)=1-\sqrt{1-2x-2y+y^2}=1-\left(1-(2x+2y-y^2) \right)^{1/2}.\]
\noindent By the Binomial theorem, we know that
\begin{eqnarray*} \left((1-(2x+2y-y^2) \right)^{1/2} & = 
& \sum_{m \geq 0} (-1)^m {1/2\choose m}\left(2x+2y-y^2\right )^m.
\end{eqnarray*}

When computing the $m$th power of $(2x+2y-y^2)$, let us consider the summand $(2x)^i(2y)^j(-y^2)^{m-i-j}$. 
The number of such summands is clearly ${m\choose i,j,m-i-j}$. Such a summand will contain $x$ and $y$ raised to the
same exponent $n$ if and only if $i=n$ and $n=j+(2m-2n-2j)$, that is, when $3n=2m-j$. In particular,
$\left(2x+2y-y^2\right )^m$ will contain a constant multiple of $x^ny^n$ if and only if 
$1.5n \leq m \leq 2n$. Therefore, if we denote the coefficient of $x^ny^n$ in $1-\sqrt{1-2x-2y+y^2} $ by
$a_{n,n}$, then routine simplifications lead to the formulas
\begin{equation}
\label{diagz} \hbox{diag}A_{(H,\phi)}(z)=\sum_{n\geq 0}a_{n,n}z^n = \sum_{n\geq 1} \left(\sum_{m= 3n/2}^{2n}
  {1/2 \choose m}{m\choose n}{m-n\choose 2m-3n} 4^{m-n}\right) z^n .\end{equation}
and
\begin{equation} a\left (G_{(H,\phi)}(n,n)\right ) =  \sum_{m= 3n/2}^{2n}
  {1/2 \choose m}{m\choose n}{m-n\choose 2m-3n} 4^{m-n}.
\end{equation}

We would like to point out that even if we have given an exact formula for $a_{n,n}$, the question of how fast
the $a_{n,n}$ grows is far from being answered. We will discuss that question in Section~\ref{growthrate}.
 Furthermore, $A_{(H,\phi)}(z)$ is algebraic, so by 
Theorem \ref{algtodfinmult}, it is $d$-finite. Therefore, Theorem \ref{diagdfin} implies that 
$ \hbox{diag}A_{(H,\phi)}(z)$ is $d$-finite. So, by Theorem \ref{dfinprec}, the coefficients of 
$\hbox{diag}A_{(H,\phi)}(z)$ must satisfy a polynomial recurrence relation. What is that relation?
We will return to this question in Section \ref{recrel}. 

\begin{example} \label{bipartite}
Let $H$ be a graph on vertex set $\{u,v\}$, having one edge, the edge $uv$, and set $\phi(u)=\phi(v)=0$.
Then, as we have seen in Corollary \ref{multipartite}, we have
\[A(x,y)=1-\sqrt{(1-x)^2+(1-y)^2-1} \]
for the generating function of the number of assembly trees of a complete bipartite graph.
\end{example}

Determining the coefficients of $\hbox{diag}A(x,y)$ is much more difficult than it was for the bivariate generating
function in Example (\ref{easyex}) since directly applying the Binomial theorem would simply lead to formulae that are too complicated
to be useful.  Nevertheless, some powerful techniques recently developed by Doron Zeilberger and Moa Apagodu will enable
us to determine these numbers. We will do this in Section \ref{recrel}. The same is true for the case of
the complete {\em tripartite} graph, which is the subject of the next example. 

\begin{example} \label{tripartite}
Let $H$ be a graph on vertex set $\{u,v,w\}$, having edge set $\{uv, uw, vw\}$, and set 
$\phi(u)=\phi(v)=\phi(w)=0$.  Then, as we have seen in Corollary \ref{multipartite}, we have
\[A(x,y,z)=1-\sqrt{(1-x)^2+(1-y)^2+(1-z)^2-2}.\]
\end{example}

\subsection{Finding Recurrence Relations} \label{recrel}
In this section, our goal is to find polynomial recurrence relations for the one-variable generating functions studied in 
Examples \ref{easyex}, \ref{bipartite} and \ref{tripartite}. As we explained in the discussion of Example \ref{easyex},
such recurrence relations exist, since our functions are diagonals of algebraic, and hence, $d$-finite power series in
several variables.

Until recently, the best available technique at this point would have been simply guessing. That is, one would have had
to assume that the sought  polynomial recurrence relation does not consist of too many terms, and does not involve polynomials of too high degrees, and then have a software package look for a suitable recurrence relation within
those limits. One major problem with this approach is that {\em even if} the software package does return a 
recurrence relation that is satisfied by all available data points, a theoretical proof that the obtained recurrence
relation is satisfied by {\em all natural numbers} $n$ is still lacking. 

The breakthrough is achieved by the following theorem of Zeilberger and Apagodu. 
(We present a simplified version of the theorem. The interested reader should consult \cite{Moa-Zeilberger}
for the full version.)

\begin{theorem} \label{maz}
Let \[F(n;x_1,x_2,\cdots ,x_d)=\prod_{p=1}^P \left(S_p(x_1,x_2,\cdots ,x_d)^{\alpha_p}\right)
\cdot  \left({s(x_1,x_2,\cdots ,x_d)}{t(x_1,x_2,\cdots,x_d)}\right)^n\]
where the $\alpha_p$ are commuting indeterminates, and 
where the $S_p$, $s$ and $t$ are elements of $\mathbb{C}[x_1,x_2,\cdots ,x_d]$.

Then there exists a non-negative integer $L$,
 there exist $L+1$ polynomials in $n$, $e_0(n),e_1(n),\cdots ,e_L(n)$, not all zero, and there exist
$d$ rational functions $R_i(n;x_1,x_2,\cdots x_d)$  (with $i=1,2,\cdots ,d$) such that
the functions
\[G_i(n;x_1,x_2,\cdots x_d):=R_i(n;x_1,x_2,\cdots x_d) F(n;x_1,x_2,\cdots,x_d) \]
satisfy the equations
\begin{equation} \label{polirek}
\sum_{i=0}^L e_i(n) F(n+i;x_1,x_2,\cdots,x_d) = \sum_{i=1}^d D_{x_i} G_i(n;x_1,x_2,\cdots x_d).
\end{equation}
Furthermore, there exists a constant $N=N(deg(s),deg(t),\sum_{p=1}^P deg(S_p))$ such that
$L\leq N$ and $deg(e_i)\leq N$, and therefore, the polynomials $e_i(n)$ can be explicitly computed.  
\end{theorem}

Note that Theorem \ref{maz} eliminates the need for "guessing" the polynomial recurrence satisfied by
the functions $F(n+i;x_1,x_2,\cdots ,x_d)$. Indeed, because of the existence of the upper bound $N$ for the number
and maximum degree of the polynomials $e_i(n)$,  finding the polynomials $e_i(n)$ in (\ref{polirek}) is simply the
question of solving a (possibly huge) system of linear equations. Furthermore, by Theorem \ref{maz} we know that
a polynomial recurrence relation (\ref{polirek}) exists, so once we have enough equations in our system to 
obtain a {\em unique} solution for the vector $e(n)=(e_0(n),e_1(n),\cdots ,e_L(n))$, we can be sure that the polynomial recurrence
relation defined by $e(n)$ is indeed a recurrence relation that is satisfied by the $F(n+i;x_1,x_2,\cdots ,x_d)$ for
{\em all} $n$. So Theorem \ref{maz} does take care of both problems we had with simply guessing a 
polynomial recurrence relation. 

Now apply  Theorem \ref{maz} to Example \ref{easyex} by letting $d=2$, $x_1=x$, and $x_2=y$. Set $S_1(x,y)=1-2x-2y+y^2$, with $\alpha_1=1/2$,
and, crucially, $S_2(x,y)=xy$ with $\alpha_2=-1$. Finally, set $s(x,y)=1$ and $t(x,y)=xy$. This leads to 
\begin{equation} \label{beforecauchy} F(n;x,y)=\frac{\sqrt{1-2x-2y+y^2}}{x^{n+1}y^{n+1}}.\end{equation}
Consider the right-hand side as a power series in two variables and integrate both sides of  (\ref{beforecauchy}) on a two-dimensional polydisk whose interior contains 0. By the two-variable residue theorem,  we see
that the only summand in the right-hand side whose integral does not vanish is $a_{n,n}\frac{1}{x}\cdot \frac{1}{y}$.
 In fact, by the residue theorem, we get
\[\int F(n,x,y) = -4\pi^2 a_{n,n}.\]

Therefore,  the polynomial recurrence relation (\ref{polirek}) is equivalent to a polynomial recurrence relation for 
the numbers $a_{n,n}$. In order to obtain this recurrence relation, we need to solve a large system of 
linear equations. Fortunately, Doron Zeilberger's software package, SMAZ \cite{ZeilRec} can do that for us. The
result is the following theorem.  

\begin{theorem} \label{anrec}
Let $a_n=a_{n,n}$ be the coefficient of $z^n$ in (\ref{diagz}). Recall that $a_n$ counts assembly trees of graphs
studied in Example \ref{easyex}.  Then we have $a_0=0$, $a_1=1$, and
\begin{equation} \label{recdiagz} a_{n+1} = \frac{3}{2} \cdot \frac{(3n-1)(3n+1)}{(n+1)^2 } a_n .\end{equation}
\end{theorem}

(Note that the above discussion computes a recurrence relation for the coefficients of 
$\rm{diag}\sqrt{1-2x-2y+y^2}$ and not
$\rm{diag}(1-\sqrt{1-2x-2y+y^2})$, but it is obvious that starting with the coefficient of $z$, the coefficients of  these two power series will satisfy the same recurrence relation.)

Similarly,  we can apply Theorem \ref{maz} to obtain a polynomial recurrence relation for the coefficients
 $b_n$ of $z^n$ in $\hbox{diag}A(z)$, where $A(x,y)$ is the bivariate generating function
of the number of assembly trees of complete bipartite graphs as computed in Example \ref{bipartite}. 
We assign the same values to the various parameters as we did immediately preceding (\ref{beforecauchy}),
except that we set $S_1(x,y)=x^2+y^2-2x-2y+1$. The result is the following. 

Let $[z^n]g(z)$ denote the coefficient of $z^n$ in the power series $g(z)$. 

\begin{theorem} \label{bnrec} Let $b_n=[z^n]\rm{diag}\left(1-\sqrt{x^2+y^2-2x-2y+1}\right)$. Note that
$b_n$ is the number of assembly trees of $K_{n,n}$, the graph studied in Example \ref{bipartite}. 
Then we have $b_0=0$, $b_1=1$, $b_2=5/2$, and 
\begin{equation} \label{recdiagcobi} 
b_{n+2}=\frac{2(6n^2+12n+5) }{(n+2)^2}b_{n+1} - \frac{n(2n-1)(2n+3)}{(n+2)^2(n+1)}b_n
\end{equation}
if $n\geq 1$. 
\end{theorem}

Finally, let $c_n=[t^n]\hbox{diag}A(t)$, where $A(x,y,z)$ is the trivariate generating function
of the number of assembly trees of complete tripartite graphs as computed in Example \ref{tripartite}.
We can then use Theorem \ref{maz} with $d=3$, $x_1=x$, $x_2=y$, $x_3=z$. Set $S_1(x,y,z)=A(x,y,z)$ as defined
in Example \ref{tripartite} with $\alpha_1=1/2$,
and  $S_2(x,y,z)=xyz$ with $\alpha_2=-1$. Finally, set $s(x,y,z)=1$, and $t(x,y,z)=xyz$. We then get the
following result. 

\begin{theorem}  Let $c_n=[t^n]\rm{diag}A(t)$, where \[A(x,y,z)=1-\sqrt{(1-x)^2+(1-y)^2+(1-z)^2-2}.\]
 Then
we have $c_0=0$, $c_1=3$,  $c_2=84$,  $c_3=4935$, and
\begin{equation} \label{recdiagcotri} c_{n+3}=r_2(n)c_{n+2}+ r_1(n)c_{n+1}+r_0(n)c_n, \end{equation}
if $n\geq 1$. Here the $r_i$ are explicitly known rational functions of $n$, with numerators and denominators of degree 11 for $r_0$, 
degree ten for $r_1$, and degree nine for $r_2$. 
\end{theorem}

\subsection{Finding Growth Rates} \label{growthrate}
As far as determining the growth rate of the sequence $a_1,a_2,\cdots $,  recurrence relation
(\ref{recdiagz})  is much more useful than the explicit formula that
we found for $a_n$  in (\ref{diagz}). Indeed, the {\em exponential} growth rate of $a_n$ is easy to read off
from (\ref{recdiagz}). It is routine to prove that 
\[ \lim_{n\rightarrow \infty} \sqrt[n]{a_n} = 13.5 .\]

Determining the growth rate of $a_n$  at a higher level of precision is much more difficult.  
The theoretical foundation of this computation is the paper \cite{Zeilberger-Wimp} by Doron Zeilberger and
Jet Wimp. In that paper, the authors consider a more general setup, but in the examples we study, their method
simplifies to the following. 

Let $t_n$ be a sequence for which a polynomial recurrence relation is known, and of which we want to compute the
asymptotics. Try to obtain $t_n$ in the form
\[t_n = E_nK_n,\]
where \[E_n= e^{\mu_0 n\ln n + \mu_1n\ln n}n^{\theta}, \]
and \[K_n=\exp \left( \alpha_1n^{\beta} + \alpha_2 n^{\beta -\frac{1}{\rho}} + \alpha_3 n^{\beta -
\frac{2}{\rho}+\cdots } \right).
\]
Here $\alpha_1\neq 0$, $\beta=j/\rho$, and $0<j<\rho$. 
This decomposition leads to the formula
\[\frac{t_{n+k}}{t_n}=n^{\mu_0k}\lambda^k \left(1+\frac{k\theta + k^2\mu_0/2}{n}+\cdots \right) 
\cdot \exp\left(\alpha_1n^{\beta}kn^{\beta-1} + \alpha_2\left(\beta-\frac{1}{\rho}k\right)n^{\beta -1 -\frac{1}{\rho}}+ \cdots
\right),\]
where $\lambda=e^{\mu_0+\mu_1}$. 
Then, using the polynomial recurrence relation for $t_n$, determine the parameters in $E_n$ and $K_n$, obtaining
this way an {\em exact} formula (in the form of an infinite sum) for $t_n$. 
This computation can certainly be long and teadious, but the software package AsyRec \cite{asyrec} can carry it out. 
For the sequence $a_n$, we obtain the following result. 

\begin{theorem} Let $a_n$ be defined as in Theorem \ref{anrec}. Note that $a_n$ is the number of assembly trees of the graph discussed in Example \ref{easyex}. Then we have
\[a_n = \frac{13.5^n}{n^2}\cdot \left(1+ \frac{1}{9n}+\frac{5}{81n^2}+\cdots \right).\]
In particular, $a_n\sim \frac{13.5^n}{n^2}$, and therefore, by Theorem \ref{jungen}, 
$\rm{diag}A_{(H,\phi)}(z)$ is {\em not} algebraic. 
\end{theorem}

Applying the same method for the polynomial recurrence relation proved for the numbers $b_n$ in 
Theorem \ref{bnrec}, we get the following asymptotic expressions. 

\begin{theorem}
Let $b_n$ be defined as in Theorem \ref{bnrec}, that is, let $b_n$ be the number of assembly trees of the complete
bipartite graph $K_{n,n}$, which was studied in Example \ref{bipartite}. Then we have 
\[b_n=\frac{(6+4\sqrt{2})^n}{n^2} \left(1+\frac{35}{8n} - \frac{5}{32}\frac{\sqrt{2}}{n} + \cdots \right) .\]
 In particular, $b_n\sim \frac{(6+4\sqrt{2})^n}{n^2}$, and therefore, by Theorem \ref{jungen}, 
$\rm{diag}A_{x,y}(z)$ is {\em not} algebraic. 
\end{theorem}

\section{Questions} \label{secOpen}  There are reasonable alternatives to the edge gluing rule for defining
an assembly tree.  Two possibilities are the following.  For each of these two rules, the problem is again to enumeratate
the number of assembly trees for interesting graphs.  \B

\begin{definition} An assembly tree for a connected graph $G$ using the  {\em connected gluing rule} is an assembly tree for $G$ (satisfying properties (1-4) in Section~\ref{secIntro}), and also satisfying the additional property:
\begin{enumerate}
\item[5.] For each node, the graph induced by the vertices in the label is connected. 
\end{enumerate} 
\end{definition} \m
\noindent This rule is  less restrictive than the edge gluling rule.  In particular, an assembly tree of a graph $G$ is not necessarily a binary tree.  At the extreme is the assembly tree of depth $1$ for which every vertex of $G$ is a child of the root.  \B

\begin{figure}[htb] \label{fig4}
\begin{center}
\includegraphics[width=9cm, keepaspectratio]{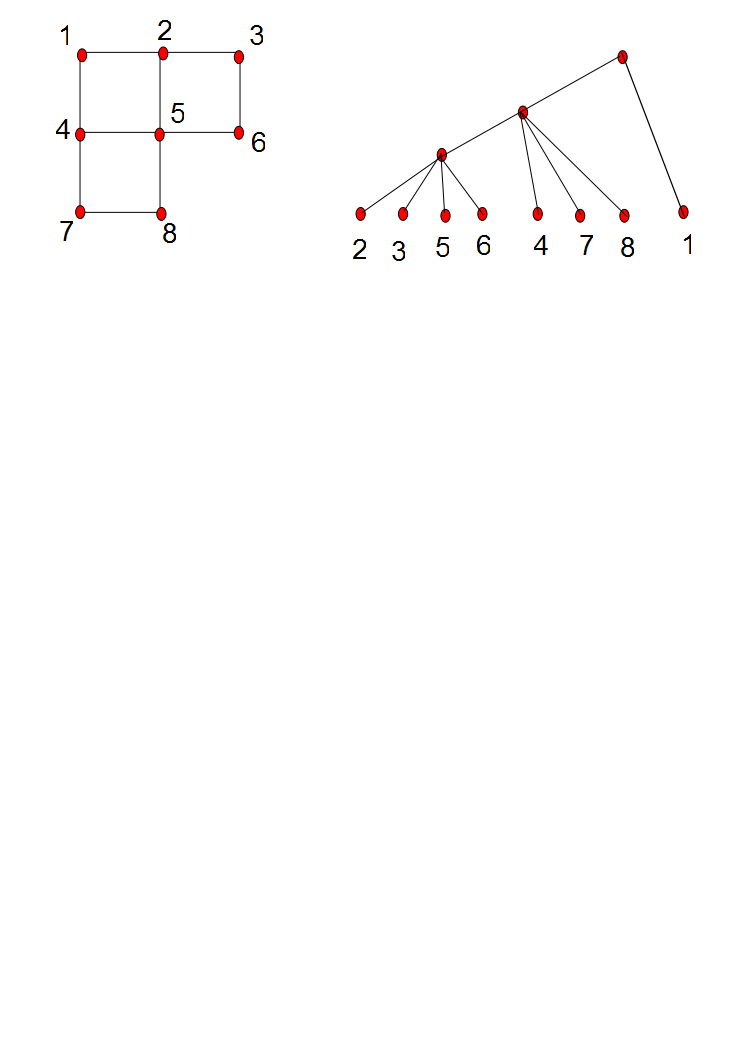}
\vskip -9cm
\caption{A plane graph and an assembly tree using the face gluing rule.}
\end{center} 
\end{figure}

\begin{definition} In this definition, we denote each face of a plane graph by the set of vertices on that face. An assembly tree for a connected plane graph $G$
using the {\em face gluing rule} is an assembly tree for $G$ satisfying the additional property. 
\begin{enumerate}
\item [5.] For each internal node $U$ there is a face $F$ such that $C\cap F \neq \emptyset$ for each
$C \in c(U)$ and $F \subseteq \bigcup \, c(U)$.   
\end{enumerate}
\end{definition} \m
\noindent Figure~\ref{fig4} shows a plane graph and one assembly tree using the face gluing rule.

\begin{figure}[htb] \label{fig5}
\begin{center}
\includegraphics[width=9cm, keepaspectratio]{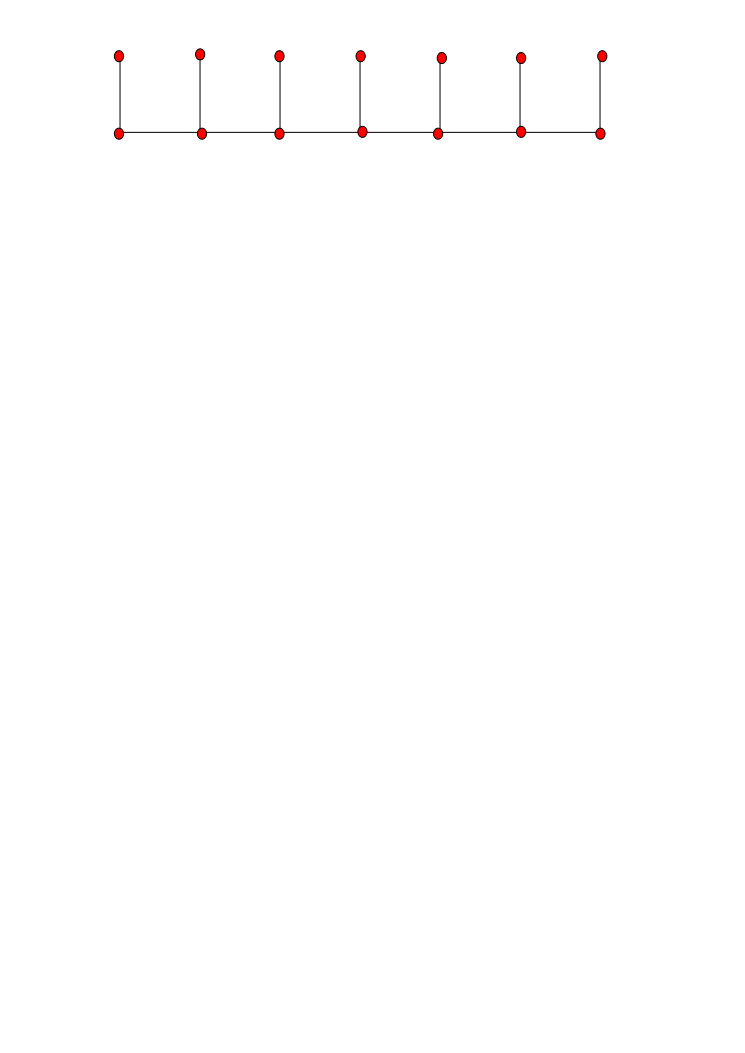}
\vskip -10.5cm
\caption{Caterpillar graph $D_7$.}
\end{center} 
\end{figure}

All the graphs in Sections~\ref{secPathCycle} and \ref{secH-partite} for which we were able to compute the number of assembly trees share a common property: the number of connected induced subgraphs, up to unlabeled isomorphism type, is small.  For example, the number of connected induced subgraphs of the path $P_n$ is $n$.  For the complete graph $K_n$,  it is $n$, and for the complete bipartite graph $K_{m,n}$ it is $mn$.  With $(H,\phi)$ fixed, the number of
connected induced subgraphs of any $(H,\phi)$-graph is again polynomial in the parameters $n_1,n_2, \dots , n_N$.   
For results on graphs with few isomorphism types of induced subgraphs  see \cite{AB} and \cite{EH}. So the question arises: is it possible to enumerate the number of assembly trees (by the edge gluing rule) for a family of graphs for which the number of isomorhism types of induced subgraphs is not small.  In particular, consider the family of caterpillar graphs $D_n$ on $2n$ vertices as shown in Figure~\ref{fig5}. The number of isomorphism types of linduced subgraphs is clearly large, exponential in $n$.   \m

{\em Question.}  Is there a reasonable enumeration of the number of assembly trees for the family $D_n$?

\vskip 1.5 cm 
\centerline{{\em Acknowledgment}}

We are indebted to Doron Zeilberger, who has  made us aware of his results in \cite{Moa-Zeilberger}, and has patiently explained
his method to us. We are also grateful to Richard Stanley, Leonard Lipshitz,  Robin Pemantle and Vincent Vatter for valuable suggestions.

\end{document}